\documentclass[12pt,a4paper,dvipdfm]{article}
\usepackage{graphicx}

\usepackage[ruled,lined,algonl,boxed]{algorithm2e}

\begin{document}
\newcommand{\up}{\vspace*{-0.05cm}}
\newcommand{\qed}{\hfill$\rule{.05in}{.1in}$\vspace{.3cm}}
\newcommand{\pf}{\noindent{\bf Proof: }}
\newtheorem{thm}{Theorem}
\newtheorem{lem}{Lemma}
\newtheorem{prop}{Proposition}
\newtheorem{prob}{Problem}
\newtheorem{quest}{Question}
\newtheorem{ex}{Example}
\newtheorem{cor}{Corollary}
\newtheorem{conj}{Conjecture}
\newtheorem{cl}{Claim}
\newtheorem{df}{Definition}
\newtheorem{rem}{Remark}
\newcommand{\beq}{\begin{equation}}
\newcommand{\eeq}{\end{equation}}
\newcommand{\<}[1]{\left\langle{#1}\right\rangle}
\newcommand{\be}{\beta}
\newcommand{\ee}{\end{enumerate}}
\newcommand{\Bul}{\mbox{$\bullet$ } }
\newcommand{\al}{\alpha}
\newcommand{\ep}{\epsilon}
\newcommand{\si}{\sigma}
\newcommand{\om}{\omega}
\newcommand{\la}{\lambda}
\newcommand{\La}{\Lambda}
\newcommand{\Ga}{\Gamma}
\newcommand{\ga}{\gamma}
\newcommand{\im}{\Rightarrow}
\newcommand{\2}{\vspace{.2cm}}
\newcommand{\es}{\emptyset}

\title{\huge\bf On General Frameworks and Threshold Functions for Multiple Domination}
\author{Vadim Zverovich\footnote{e-mail: {\tt vadim.zverovich@uwe.ac.uk}}\\
{\footnotesize University of the West of England, Bristol, BS16
1QY, UK}
\\
}
%\date{arXiv:}
%\date{}
\maketitle

\begin{abstract}

We consider two general frameworks for multiple domination, which are called 
$\left\langle \mathbf{r},\mathbf{s} \right\rangle$-domination and parametric domination.
They generalise and unify $\{k\}$-domi\-nation, $k$-domination, total $k$-domination and $k$-tuple domination.
In this paper, known upper bounds for the classical domination are generalised 
for the $\left\langle \mathbf{r},\mathbf{s} \right\rangle$-domina\-tion and parametric domination numbers.
These generalisations imply new upper bounds for the $\{k\}$-domination and total $k$-domination numbers. 
Also, we study threshold functions, which impose additional restrictions on the minimum vertex degree,
and present new upper bounds for the aforementioned numbers. Those bounds extend similar known results for $k$-tuple domination and total $k$-domination.

\bigskip
\noindent {\it Keywords:} $\left\langle \mathbf{r},\mathbf{s} \right\rangle$-domination, parametric domination, 
$\{k\}$-domination, $k$-domination, total $k$-domination, $k$-tuple domination, upper bounds, threshold functions 
\end{abstract}

\section{Introduction}
\label{intro}

All graphs will be finite and undirected without loops and
multiple edges. If $G$ is a graph of order $n$, then
$V(G)=\{v_1,v_2,...,v_n\}$ is the set of vertices in $G$ and $d_i$
denotes the degree of $v_i$. Let $N(x)$ denote the neighbourhood of a
vertex $x$. Also let $ N(X)=\cup_{x\in X} {N(x)} $ and $
N[X]=N(X)\cup X.$ Denote by $\delta=\delta(G)$ and $\Delta=\Delta(G)$ the
minimum and maximum degrees of vertices of $G$, respectively. 
The following well-known definitions can be found in \cite{hay1}.
A set $X$ is called a {\it dominating set} if every vertex not in
$X$ is adjacent to a vertex in $X$. The minimum cardinality of a
dominating set of $G$ is the {\it domination number} $\gamma(G)$.
A set $X$ is called a {\it $k$-dominating set} if every vertex not
in $X$ has at least $k$ neighbours in $X$. 
The minimum cardinality
of a $k$-dominating set of $G$ is the {\it $k$-domination number}
$\gamma_k(G)$. 
A set $X$ is called a {\it $k$-tuple dominating set} of $G$ if for every
vertex $v\in V(G)$, $|N[v]\cap X|\ge k$. The minimum cardinality
of a $k$-tuple dominating set of $G$ is the {\it $k$-tuple
domination number} $\gamma_{\times k}(G)$. The $k$-tuple
domination number is only
defined for graphs with $\delta\ge k-1$. The 2-tuple domination
number $\ga_{\times2}(G)$ is called the {\it double domination
number} and the 3-tuple domination number $\ga_{\times3}(G)$ is
called the {\it triple domination number}. A set $X$ is called a
{\it total $k$-dominating set} of $G$ if for every vertex $v\in
V(G)$, $|N(v)\cap X|\ge k$. The minimum cardinality of
a total $k$-dominating set of $G$ is the {\it total $k$-domination
number} $\gamma_k^{\mathrm t}(G)$. The total $k$-domination number
is only defined for graphs with $\delta\ge k$. Note that
$\gamma_1^{\mathrm t}(G)$ is the well-known total domination
number $\gamma_{\mathrm t}(G)$.
A function $f: V(G) \rightarrow \{0,1,...,k\}$ is called $\{k\}$-{\em dominating}  if 
$f(N[v_i]) \ge k$ for all $i=1,...,n$.
The $\{k\}$-{\em domination number} of a graph $G$, 
denoted by $\gamma_{\{k\}}(G)$,
is the smallest weight of a $\{k\}$-dominating function of $G$.

Various problems in ad hoc networks, biological networks,
distributed computing and social networks can be modelled by dominating sets in graphs.
For example, for balancing efficiency and fault
tolerance in wireless sensor networks multiple dominating sets can be used \cite{dai1}.

The following fundamental result was independently proved by Alon
and Spencer \cite{alo1}, Arnautov \cite{arn1}, Lov\'{a}sz
\cite{lov1} and Payan \cite{pay1}. Notice that a simple
deterministic algorithm to construct a dominating set satisfying
bound (\ref{classical}) can be found in \cite{alo1}.

\begin{thm} [\cite{alo1,arn1,lov1,pay1}] \label{arn}
For any graph $G$,
\begin{equation}\label{classical}
\gamma(G) \le {\ln(\delta+1)+1 \over \delta+1} n.
\end{equation}
\end{thm}

Alon \cite{alo2} proved that the bound of Theorem \ref{arn} is asymptotically best possible. 
The bound (\ref{Car-Rod}) in Theorem \ref{Caro} is asymptotically same as bound (\ref{classical}), 
even though  (\ref{Car-Rod}) is stronger for small values of $\delta$. 

\begin{thm}[\cite{hay1}, p.\,48]
\label{Caro} For any graph $G$ with $\delta\ge 1$,
\begin{equation}\label{Car-Rod}
\gamma(G) \le \left(1-{\delta \over
{(1+\delta)^{1+1/\delta}}}\right) n.
\end{equation}
\end{thm}

In this paper, we consider Cockayne's and Favaron's general frameworks for 
$\left\langle \mathbf{r},\mathbf{s} \right\rangle$-domination and parametric domination, respectively.
They generalise and unify $\{k\}$-domination, $k$-domination, total $k$-domination and $k$-tuple domination.
In the next two sections, the classical upper bounds (\ref{classical}) and (\ref{Car-Rod}) are generalised 
for the $\left\langle \mathbf{r},\mathbf{s} \right\rangle$-domination and parametric domination numbers.
These generalisations imply new upper bounds for the $\{k\}$-domination and total $k$-domination numbers (see Section 4).

In Section 5, we study threshold functions, which impose additional restrictions on the minimum vertex degree,
and present new upper bounds for the above numbers. Those bounds extend similar known results for $k$-tuple domination and total $k$-domination.

Note that the probabilistic constructions used in the proofs of the theorems on 
$\left\langle \mathbf{r},\mathbf{s} \right\rangle$-domination
imply randomized algorithms for finding an $\mathbf{s}$-dominating $\mathbf{r}$-functions, 
whose weights satisfy the bounds of the corresponding theorems with positive probability.
A similar statement is true for the theorems devoted to parametric domination.

\section{\textsc{Cockayne's} framework for $\left\langle \mathbf{r},\mathbf{s} \right\rangle$-domination}

Cockayne introduced in \cite{coc1} an interesting framework for domination in graphs.
Let $V(G)=\{v_1,...,v_n\}$ denote the vertex set of a graph $G$, and let 
$\mathbf{r}=(r_1,...,r_n)$ and $\mathbf{s}=(s_1,...,s_n)$ be $n$-tuples of nonnegative integers, i.e. 
$r_i\in \textbf{N}$ and $s_i\in \textbf{N}$.
A function $f: V(G) \rightarrow\textbf{N}$ is called an $\mathbf{r}$-{\em function} of $G$ if
$$
f(v_i) \le r_i \quad \mbox{for all} \quad i=1,...,n.
$$
Let 
$$
f[v_i] = \sum_{u\in N[v_i]} f(u).
$$
An $\mathbf{r}$-function $f$ is $\mathbf{s}$-{\em dominating}  if 
$$
f[v_i] \ge s_i \quad \mbox{for all} \quad i=1,...,n.
$$
The {\em weight} of a function $f$ is denoted by $|f|$ and defined by
$$
|f| = \sum_{i=1}^n f(v_i). 
$$

The $\left\langle \mathbf{r},\mathbf{s} \right\rangle$-{\em domination number} of a graph $G$, 
denoted by $\gamma\left\langle \mathbf{r},\mathbf{s} \right\rangle (G)$,
is the smallest weight of an $\mathbf{s}$-dominating $\mathbf{r}$-function of $G$.
As pointed out in \cite{coc1}, such functions exist if and only if
$$
\sum_{v_j\in N[v_i]} r_j \ge s_i \quad \mbox{for all} \quad i=1,...,n.
$$
It is not difficult to see that $\left\langle \mathbf{r},\mathbf{s} \right\rangle$-domination unifies and generalises the classical domination, $k$-tuple domination and $\{k\}$-domination 
if we put $r_i=s_i=1$;  $r_i=1, s_i=k$; and $r_i=s_i=k$ for all $i=1,...,n$, respectively.

Let us denote
$$
\tau = \min\{r_1,...,r_n\}, \quad s = \max\{s_1,...,s_n\}, \quad r = \left\lfloor {s \over \delta+1} \right\rfloor +1,
$$
%In what follows, we assume that $(\delta+1)r \ge s$ to guarantee the existence of an %$\mathbf{s}$-dominating $\mathbf{r}$-function. Note that any $(s,...,s)$-dominating 
%$(r,...,r)$-function is also an $\mathbf{s}$-dominating $\mathbf{r}$-function. 
%Hence 
%$$
%\gamma\left\langle \mathbf{r},\mathbf{s} \right\rangle (G) \le
%\gamma\left\langle (r,...,r),(s,...,s) \right\rangle (G).
%$$
%Thus, for simplicity of presentation, we formulate our result for the case when 
%$\mathbf{r}=(r,...,r)$ and $\mathbf{s}=(s,...,s)$ are $n$-tuples of nonnegative integers numbers $r$ and %$s$. 
$$
\theta = (\delta+1)r - s 
 \quad \mbox{and} \quad
 {B}_t =  \pmatrix{(\delta+1)r \cr t}.
$$

The following theorem provides an upper bound for the $\left\langle \mathbf{r},\mathbf{s} \right\rangle$-domination number of a graph. 

\begin{thm} \label{ub1}
For any graph $G$ of order $n$ with $r\le \tau$ and $\rho=1/\theta$,
\begin{equation}
\gamma\left\langle \mathbf{r},\mathbf{s} \right\rangle (G) \le \left(1-{ (r \rho)^\rho 
\over (1+\rho)^{1+\rho}\; {B}_{s-1}^{\rho} }\right) rn.
\end{equation}
\end{thm}

\begin{pf}
For each vertex $v\in V(G)$, we select $\delta$ vertices from
$N(v)$ and denote the resulting set together with the vertex $v$ by $N'[v]$. 
Thus, $|N'[v]|=\delta+1$.
For $i=1,2,...,r,\;$ let $a_i(v)$ be a (0,1)-function on the set $V(G)$ such that it assigns ``1'' to every vertex of $G$ independently with probability 
$$
p =
1-\left({r\over (1+\theta)B_{s-1}}\right)^{1/\theta}.
$$
Let us define an $\mathbf{r}$-function $a(v)$ as follows:
$$
a(v)=\sum_{i=1}^r a_i(v).
$$

%For $i=1,2,...,r$, let $A_i$ be a set formed by an independent choice of vertices
%of $G$, where each vertex is selected with the probability 
%$$
%p = 1-\left({r\over (1+\theta)B_{s-1}}\right)^{1/\theta}.
%$$
%Let us form a multiset $A=\cup_{i=1}^r$, and let $a(v)$ be a function on the set $V(G)$ such that it is %equal to the multiplicity of the vertex $v$ in $A$. Thus, $0\le a(v)\le r$ for any vertex $v$.

For $m=0,1,...,s-1$, we denote 
$$ 
C_m = \{v\in V(G) : \sum_{u\in N'[v]} a(u)=m\}.
$$ 

\begin{cl} \label{claim1}
For each set $C_m$, there exists a function $c_m: V(G) \rightarrow\textbf{N}$ such that
\begin{equation}
\label{e1}
|c_m|\le (s-m)|C_m|
\end{equation}
and for any vertex $v\in C_m$,
\begin{equation}
\label{e2}
a(v)+c_m(v) \le r   \quad \mbox{and} \quad \sum_{u\in N'[v]} c_m(u) \ge s-m.
\end{equation}
\end{cl}

\begin{pf}
Let us initially put $c_m(v)=0$ for all $v\in V(G)$.
Then, for each vertex $v\in C_m$, we redefine $c_m$ in the set $N'[v]$ as follows:

{\it Case 1:} Suppose that $c_m(u)=0$ for any vertex $u\in N'[v]$.
Note that because $\sum_{u\in N'[v]} a(u)=m$, the ``spare capacity'' in $N'[v]$ is 
$(\delta +1)r-m > s-m$, i.e. the weight of $c_m$ can be increased in $N'[v]$ by $s-m$ units.
Thus, we can obviously redefine $c_m$ in $N'[v]$ in such a way 
that 
\begin{equation}
\label{e3}
\sum_{u\in N'[v]} c_m(u)=s-m
\end{equation}
and 
\begin{equation}
\label{e4}
a(u)+c_m(u) \le r \quad \mbox{for any} \quad u\in N'[v].
\end{equation}
In this case, we increased the weight of $c_m$ in $N'[v]$ by $s-m$ units.

{\it Case 2:} Assume that $c_m(u)>0$ for some $u\in N'[v]$, but 
$$\sum_{u\in N'[v]} c_m(u) = \psi < s-m,$$
where $\psi \ge 1$.
In this case, we can increase the weight of $c_m$ in $N'[v]$ by $s-m-\psi$ units
to make sure that (\ref{e3}) and  (\ref{e4}) hold.

{\it Case 3:} Suppose now that $c_m(u)>0$ for some $u\in N'[v]$ and 
$$\sum_{u\in N'[v]} c_m(u) \ge s-m.$$
In this case, we do not change the weight of $c_m$ in $N'[v]$.

Thus, when constructing the function $c_m$, we increased its weight at most $|C_m|$ times by at most $s-m$ units, and so (\ref{e1}) is true. The inequalities (\ref{e2}) are also  true by construction.

\qed
\end{pf}

Let us define the function $f$ on the set $V(G)$ as follows: 
$$ 
f(v) = a(v) +  \max _{0\le m\le s-1} c_m(v).
$$
By Claim \ref{claim1}, the function $f$ is an $(s,...,s)$-dominating 
$(r,...,r)$-function. Hence, it is also an $\mathbf{s}$-dominating $\mathbf{r}$-function because
$s\ge s_i$ and $r\le\tau\le r_i$ for all $i=1,...,n$.
Also, 
$$ 
f(v) \le a(v) +  \sum_{m=0}^{s-1}  c_m(v).
$$
The expectation of $|f|$ is as follows:
\begin{eqnarray*}
\mathbf{E}[|f|] &\le& \mathbf{E}\left[|a| + \sum_{m=0}^{s-1}|c_m| \right]
\le \sum_{i=1}^r\mathbf{E}[|a_i|] + \sum_{m=0}^{s-1}(s-m)\mathbf{E}[|C_m|].
\end{eqnarray*}
We have
\begin{eqnarray*}
\mathbf{E}[|C_m|]  =  \sum_{v\in V(G)} \mathbf{P}[v\in C_m]
   =  \sum_{i=1}^n  p^m (1-p)^{(\delta+1)r-m} \pmatrix{(\delta+1)r \cr m}
 = p^{m} (1-p)^{(\delta+1)r-m} {B}_m n.
\end{eqnarray*}
Thus,
\begin{eqnarray*}
\mathbf{E}[|f|] &\le& pnr + \sum_{m=0}^{s-1}(s-m) p^{m}(1-p)^{(\delta+1)r-m} {B}_m n \\
&=& pnr + (1-p)^{(\delta+1)r-s+1} n \sum_{m=0}^{s-1}(s-m) p^m
(1-p)^{s-m-1} {B}_{m}.\\
\end{eqnarray*}
Furthermore, for $0 \le m \le s-1$,
\begin{small}
\begin{eqnarray*}
{(s-m)}{B}_{m}  = (s-m) {(\delta+1)r \choose m} \le
\prod_{j=1}^{\theta+1} {(s-m+j-1) \over j} {(\delta+1)r \choose m}\\
 = {(\delta+1)r - m  \choose {\theta + 1}} {(\delta+1)r \choose m}
 = {s- 1 \choose m} {(\delta+1)r \choose {s-1}}  = {s- 1 \choose m}
{B}_{s-1}.
\end{eqnarray*}
\end{small}
We obtain
\begin{eqnarray*}
\gamma\left\langle \mathbf{r},\mathbf{s} \right\rangle (G) 
\le \mathbf{E}[|f|] &\le& pnr + {(1-p)}^{\theta +1}n {B}_{s-1} 
\sum_{m=0}^{s-1} {s-1 \choose m} p^m {(1-p)}^{s-m-1}\\
 &=& pnr + {(1-p)}^{\theta+1} n {B}_{s-1}\\
% &=& \left(1-{\theta\;  r^{1/\theta} 
%\over (1+\theta)^{1+1/\theta}\; {B}_{s-1}^{1/\theta} }\right) rn\\
&=& \left(1-{ (r \rho)^\rho 
\over (1+\rho)^{1+\rho}\; {B}_{s-1}^{\rho} }\right) rn,
\end{eqnarray*}
as required. The proof of the theorem is complete.\qed
\end{pf}

The proof of Theorem \ref{ub1} implies a weaker upper bound for the $\left\langle \mathbf{r},\mathbf{s} \right\rangle$-domination number. This result generalises the classical bound in Theorem \ref{arn}.

\begin{cor} \label{cor1}
For any graph $G$ of order $n$ with $r \le \tau$,
\begin{equation}
\gamma\left\langle \mathbf{r},\mathbf{s} \right\rangle (G) \le 
{\ln(\theta+1) + \ln B_{s-1} -\ln r +1 \over \theta +1} rn.
\end{equation}
\end{cor}

\begin{pf}
The proof easily follows if we use the inequality $1-p \le e^{-p}$, and then minimise the following upper bound:
$$
\gamma\left\langle \mathbf{r},\mathbf{s} \right\rangle (G) 
\le  pnr + e^{-p(\theta+1)} n {B}_{s-1}.
$$
\qed
\end{pf}

It may be pointed out that the bound of Corollary \ref{cor1} can be optimised with respect to $r$, where $r$ is now any integer between $s/(\delta+1)$ and $\tau$:
$$
\gamma\left\langle \mathbf{r},\mathbf{s} \right\rangle (G) \le 
\min _{s/(\delta+1)\le r\le \tau}\left\{{\ln(\theta+1) + \ln B_{s-1} -\ln r +1 \over \theta +1} rn\right\}.
$$

While $\left\langle \mathbf{r},\mathbf{s} \right\rangle$-domination generalises the classical domination, $k$-tuple domination and $\{k\}$-domination, the following definition generalises total domination by considering open neighbourhoods. 

An $\mathbf{r}$-function $f$ is called {\em total $\mathbf{s}$-dominating}  if 
$$
\sum_{u\in N(v_i)} f(u) \ge s_i \quad \mbox{for all} \quad i=1,...,n.
$$
The total $\left\langle \mathbf{r},\mathbf{s} \right\rangle$-{\em domination number}  
$\gamma^{\mathrm t}\left\langle \mathbf{r},\mathbf{s} \right\rangle (G)$ of a graph $G$ 
is the smallest weight of a total $\mathbf{s}$-dominating $\mathbf{r}$-function of $G$.

Let
$$
\tilde{r} = \left\lfloor {s \over \delta}\right\rfloor +1, \quad
\tilde{\theta} = \delta \tilde{r} - s, 
 \quad 
\tilde{B}_{s-1} =  \pmatrix{\delta \tilde{r} \cr s-1} \quad \mbox{and} \quad \tilde{\rho}=1/\tilde{\theta}.
$$

\begin{thm} \label{thm2}
For any graph $G$ of order $n$ with $\tilde{r} \le \tau$ and $\delta>0$,
$$
\gamma^{\mathrm t}\left\langle \mathbf{r},\mathbf{s} \right\rangle (G) \le \left(1-{ (\tilde{r} \tilde{\rho})^{\tilde{\rho}} 
\over (1+\tilde{\rho})^{1+\tilde{\rho}}\; \tilde{B}_{s-1}^{\tilde{\rho}} }\right) \tilde{r}n.
$$

Also, 
$$
\gamma^{\mathrm t}\left\langle \mathbf{r},\mathbf{s} \right\rangle (G) \le 
{\ln(\tilde{\theta}+1) + \ln \tilde{B}_{s-1} -\ln \tilde{r} +1 \over \tilde{\theta} +1} \tilde{r}n.
$$
\end{thm}

\begin{pf}
For each vertex $v\in V(G)$, we select $\delta$ vertices from
$N(v)$ and denote the resulting set $N'(v)$. 
Thus, $|N'(v)|=\delta$.
The proof now follows immediately from the proofs of Theorem \ref{ub1} and Corollary \ref{cor1} if we replace 
$N'[v]$ by $N'(v)$. 
\qed
\end{pf}

\section{Favaron's framework for parametric domination}
\label{par}

While Cockayne's framework is based on functions with prescribed properties, 
the focus of the generalisation considered in this section is on properties of
vertex sets called $(k,l)$-dominating sets. These two frameworks complement each other 
because the former does not generalise $k$-domination, while the latter does not include $\{k\}$-domination. 

The following definition with minor adaptations is due to Favaron et al \cite{fav1}.
For integers $k\ge1$ and $l\ge1$, a set $D$ is called a
$(k,l)$-{\em dominating set} of $G$ if for every vertex $v\not\in
D$,
$$|N[v]\cap D|\ge k,$$ and for every vertex $v\in D$, $$|N[v]\cap
D|\ge l.$$ The minimum cardinality of a $(k,l)$-dominating set of
$G$ is the {\em parametric domination number} $\gamma_{k,l}(G)$.
It is natural to consider the parametric domination number for graphs with 
$\delta\ge \max\{k,l-1\}$. Since $V(G)$ is a $(k,l)$-dominating set of $G$,
the parametric domination is well defined. It is easy to see that
$\gamma_{1,1}(G)$ is the domination number $\gamma(G)$,
$\gamma_{2,1}(G)$ is the 2-domination number $\gamma_{2}(G)$,
$\gamma_{2,2}(G)$ is the double domination number
$\ga_{\times2}(G)$ and
$\gamma_{1,2}(G)$ is the total domination number $\ga_{\mathrm
t}(G)$.

More generally, the parametric domination number unifies the following:

\bigskip
\begin{tabular}{|l|l|}
  \hline
  % after \\: \hline or \cline{col1-col2} \cline{col3-col4} ...
  $l=1$ & $\gamma_{k,1}(G)$ is the $k$-domination number $\gamma_k(G)$
  \\
  \hline
  $l=k$ & $\gamma_{k,k}(G)$ is the $k$-tuple domination number
$\gamma_{\times k}(G)$ \\
\hline
  $l=k+1$ & $\gamma_{k,k+1}(G)$ is the total $k$-domination number
$\gamma_k^{\mathrm t}(G)$ \\
  \hline
\end{tabular}
\bigskip

Let
$$
\varphi = \max\{k,l-1\} \quad \mbox{and} \quad b_{t}=\pmatrix{\delta \cr t}.
$$

\begin{thm}
\label{main} For any graph $G$ with $\bar\delta  = \delta - \varphi >0$,
$$
\gamma_{k,l}(G) \le 
\left(1 - {{\bar\delta} \over
({1+\bar\delta})^{1+1/{\bar\delta}} \; b_{\varphi-1}^{1/{\bar\delta}}}\right) n.
$$
\end{thm}

\pf 
For each vertex $v\in V(G)$, we select $\delta$ vertices from
$N(v)$ and denote the resulting set by $N'(v)$.
Let $A$ be a set formed by an independent choice of vertices
of $G$, where each vertex is selected with probability 
$$
p =  1-\left({1 \over (1+\bar\delta)b_{\varphi-1}}\right)^{1/\bar\delta}.
$$
For $m=0,1,...,k-1$, let us denote
$
B_m = \{v\in V(G)-A : |N'(v)\cap A|=m\}.
$
Also, for $m=0,1,...,l-2$, we denote
$
A_m = \{v\in A : |N'(v)\cap A|=m\}.
$
For each set $A_m$, we form a set $A'_m$ in the following way. For
every vertex $v\in A_m$, we take $l-m-1$ neighbours from $N'(v)-A$ and add them to $A'_m$. Such neighbours always exist because
$\delta \ge l-1$. It is obvious that $|A'_m|\le (l-m-1)|A_m|$. For
each set $B_m$, we form a set $B'_m$ by taking $k-m$ neighbours
from $N'(v)-A$ for every vertex $v\in B_m$.  Such neighbours always
exist because $\delta \ge k$. We have $|B'_m|\le (k-m)|B_m|$.

Let us construct the set $D$ as follows:
$
D = A \cup \left(\bigcup_{m=0}^{l-2} A'_m\right) \cup
\left(\bigcup_{m=0}^{k-1} B'_m\right).
$
The set $D$ is a $(k,l)$-dominating set. Indeed, if there is a
vertex $v$ which is not $(k,l)$-dominated by $D$, then $v$ is not
$(k,l)$-dominated by $A$. Therefore, $v$ would belong to $A_m$ or
$B_m$ for some $m$, but all such vertices are $(k,l)$-dominated by
the set $D$ by construction.

The expectation of $|D|$ is
\begin{eqnarray*}
\mathbf{E}[|D|] &\le& \mathbf{E}\Big[|A| + \sum_{m=0}^{l-2}|A'_m| +
\sum_{m=0}^{k-1}|B'_m|\Big]\\
&\le& \mathbf{E}[|A|] + \sum_{m=0}^{l-2}(l-m-1)\mathbf{E}[|A_m|] +
\sum_{m=0}^{k-1}(k-m)\mathbf{E}[|B_m|].
\end{eqnarray*}
We have
$
\mathbf{E}[|A|] = \sum_{i=1}^n \mathbf{P}[v_i\in A] = pn.
$
Also,
\begin{eqnarray*}
\mathbf{E}[|A_m|]  =  \sum_{i=1}^n \mathbf{P}[v_i\in A_m]  =  \sum_{i=1}^n p
\pmatrix{\delta \cr m} p^m (1-p)^{\delta-m}  =  p^{m+1} (1-p)^{\delta-m} b_m n
\end{eqnarray*}
and
\begin{eqnarray*}
\mathbf{E}[|B_m|]  =  \sum_{i=1}^n \mathbf{P}[v_i\in B_m]  =  \sum_{i=1}^n (1-p)
\pmatrix{\delta \cr m} p^m (1-p)^{\delta-m} =  p^m (1-p)^{\delta-m+1} b_m n.
\end{eqnarray*}
We obtain
\begin{eqnarray*}
\mathbf{E}[|D|] &\le& pn + \sum_{m=0}^{l-2}(l-m-1) p^{m+1} (1-p)^{\delta-m}
b_m n  + \sum_{m=0}^{k-1}(k-m)p^m (1-p)^{\delta-m+1}
b_m n\\
&\le& pn + \sum_{m=0}^{\varphi-1}(\varphi-m) p^{m+1} (1-p)^{\delta-m}
b_m n + \sum_{m=0}^{\varphi-1}(\varphi-m)p^m
(1-p)^{\delta-m+1}b_m n\\
&=& pn + \sum_{m=0}^{\varphi-1}(\varphi-m)b_m n p^{m}
(1-p)^{\delta-m}.
\end{eqnarray*}
Furthermore, for $0 \le m \le \varphi-1$,
\begin{eqnarray*}
{(\varphi-m)} b_{m} = (\varphi-m) {\delta \choose m}
 \le (\varphi-m) {\delta \choose m}
\prod_{j=2}^{\varphi-m} {(\delta-\varphi+j) \over j} \qquad\qquad\qquad\qquad\\
 = { \delta!\over m!(\varphi-m-1)!(\delta-\varphi+1)!}  
 = {\varphi- 1 \choose m} {\delta \choose {\varphi-1}}  
 = {\varphi - 1 \choose m} b_{\varphi-1}.
\end{eqnarray*}
Therefore, 
\begin{eqnarray*}
\mathbf{E}[|D|] &\le& pn + n b_{\varphi-1}(1-p)^{\delta-\varphi+1}\sum_{m=0}^{\varphi-1}{\varphi - 1
\choose m} p^{m}(1-p)^{\varphi-1-m}\\
&=& pn + n b_{\varphi-1} (1-p)^{\delta-\varphi+1}\\
&=& \left(1 - {{\bar\delta} \over
({1+\bar\delta})^{1+1/{\bar\delta}} \; b_{\varphi-1}^{1/{\bar\delta}}}\right) n.
\end{eqnarray*}
Since the expectation is an average value, there exists a
particular $(k,l)$-dominating set of the above order, as required.
The proof of Theorem \ref{main} is complete. \qed

\begin{cor} \label{kl-cor}
For any graph $G$ with $\delta \ge \varphi$,
$$
\gamma_{k,l}(G) \le  {\ln(\delta -\varphi + 1) + \ln b_{\varphi-1} +1\over \delta -\varphi + 1} n.
$$
\end{cor}

\begin{pf}
Using the inequality $1-p \le e^{-p}$, we obtain 
\begin{eqnarray*}
\mathbf{E}[|D|] \le pn + n b_{\varphi-1} e^{-p(\delta-\varphi+1)}.
\end{eqnarray*}
The proof easily follows if we minimise the right-hand side in the above inequality.
\qed
\end{pf}

The next result is similar to Theorem \ref{main} and Corollary \ref{kl-cor}, which provide 
better bounds if $l \ge k+1$. However, for small values of $l$, the bounds of Theorem \ref{main2}
are better. In what follows, we put $b_{-1}=0$.

\begin{thm} \label{main2}
For any graph $G$ with ${\widehat\delta}=\delta-\max\{k,l\}+1>0$,
$$
\gamma_{k,l}(G) \le \left(1 - {{\widehat\delta} \over
({1+\widehat\delta})^{1+1/{\widehat\delta}} \; (b_{k-1} +
b_{l-2})^{1/{\widehat\delta}}}\right) n.
$$
\2

Also, for any graph $G$ with $\delta\ge \max\{k,l-1\}$,
$$
\gamma_{k,l}(G) \le  {\ln(\widehat\delta +1) + \ln (b_{k-1} + b_{l-2}) +1\over \widehat\delta + 1} n.
$$
\end{thm}

\pf Using the same construction as in the proof of Theorem
\ref{main}, we obtain
\begin{eqnarray*}
\mathbf{E}[|D|] &\le& pn + \sum_{m=0}^{l-2}(l-m-1) p^{m+1} (1-p)^{\delta-m}
b_m n  + \sum_{m=0}^{k-1}(k-m)p^m (1-p)^{\delta-m+1}
b_m n\\
&=& pn + \sum_{m=1}^{l-1}(l-m) p^m (1-p)^{\delta-m+1} b_{m-1} n  + \sum_{m=0}^{k-1}(k-m)p^m (1-p)^{\delta-m+1}
b_m n\\
&=& pn + (1-p)^{\delta-l+2} n \Theta_1 + (1-p)^{\delta-k+2} n \Theta_2,
\end{eqnarray*}
where $$
\Theta_1 = \sum_{m=1}^{l-1}(l-m) p^m
(1-p)^{l-m-1}b_{m-1}, \quad
 \Theta_2 = \sum_{m=0}^{k-1}(k-m) p^m
(1-p)^{k-m-1} b_{m}.
$$
Now, using an approach similar to that in the proof of Theorem
\ref{main}, we can prove that
$$(l-m)b_{m-1} \le {l-1
\choose m} b_{l-2} \quad \mbox{and} \quad
(k-m)b_{m} \le {k-1 \choose m} b_{k-1}.
$$
Therefore, $\Theta_1 \le b_{l-2}$ and $\Theta_2 \le b_{k-1}$. 
We have
\begin{eqnarray*}
\mathbf{E}[|D|] &\le& pn + (1-p)^{\delta-l+2}n b_{l-2} + 
(1-p)^{\delta-k+2}n b_{k-1}\\
&\le& pn + (1-p)^{\widehat\delta+1}n (b_{l-2}+b_{k-1}).
\end{eqnarray*}
The first upper bound of the theorem is obtained by minimising the above function, 
while the second by minimising the following function: 
$
\mathbf{E}[|D|] \le  pn + e^{-p(\widehat\delta+1)}n (b_{l-2}+b_{k-1}).
$
\qed

%\vspace*{-0.8cm}

\section{Corollaries for $\{k\}$-domination, $k$-domination, total $k$-domination and $k$-tuple domination   }

$\{k\}$-Domination is a particular case of $\left\langle \mathbf{r},\mathbf{s} \right\rangle$-domination 
when $r_i = s_i = k$ for all $i=1,...,n$. Theorem \ref{ub1} and Corollary \ref{cor1} imply new upper bounds for the $\{k\}$-domination number. We have $\tau=s=k$, and hence
$$
r = \left\lfloor {k \over \delta+1} \right\rfloor +1, \quad
\theta = (\delta+1)r - k 
 \quad \mbox{and} \quad
 {B}_{k-1} =  \pmatrix{(\delta+1)r \cr k-1}.
$$

\begin{cor} \label{k-dom}
For any graph $G$ with $\delta > 0$ and $\rho=1/\theta$,
$$
\gamma_{\{k\}} (G) \le  \left(1-{ (r \rho)^\rho 
\over (1+\rho)^{1+\rho}\; {B}_{k-1}^{\rho} }\right) rn.
$$

Also, 
$$
\gamma_{\{k\}} (G) \le {\ln(\theta+1) + \ln B_{k-1} -\ln r +1 \over \theta +1} rn.
$$
\end{cor}

Similar to Corollary \ref{cor1}, the latter bound in Corollary \ref{k-dom} is weaker than the former bound, but it has a simpler formula and can be further optimised with respect to $r$ for integers between $k/(\delta+1)$ and $k$.

The special case $l=1$ in parametric domination is the $k$-domination number $\gamma_k(G)$. 
By Theorem \ref{main2}, the following result is obtained:

\begin{cor} [Gagarin et al. \cite{gag1}]
For any graph $G$ with ${\widehat\delta}=\delta-k+1$ and $\delta\ge k$,
$$
\gamma_{k}(G) \le \left(1 - {{\widehat\delta} \over
({1+\widehat\delta})^{1+1/{\widehat\delta}} \; b_{k-1}^{1/{\widehat\delta}}}\right) n
$$
and
$$
\gamma_k(G) \le  {\ln(\delta-k+2) + \ln b_{k-1} +1 \over \delta-k+2} n.
$$
\end{cor}

The following upper bounds for the total $k$-domination number 
are easily obtained if we put $l=k+1$ in Theorem \ref{main} and Corollary \ref{kl-cor}.
Notice that they also follow from Theorem \ref{thm2} if we put $r_i=1$ and $s_i=k$ for all $i=1,...,n$, 
in which case $\tilde{r}=1$, $k < \delta$, $\tilde{\theta}=\delta-k$ and $\tilde{B}_{k-1}=b_{k-1}$.

\begin{cor}
For any graph $G$ with $\bar\delta  = \delta - k >0$,
$$
\gamma_k^{\mathrm t}(G) \le 
\left(1 - {{\bar\delta} \over
({1+\bar\delta})^{1+1/{\bar\delta}} \; b_{k-1}^{1/{\bar\delta}}}\right) n.
$$

Also, for $\delta\ge k$,
$$
\gamma_k^{\mathrm t}(G) \le  {\ln(\delta -k +1) + \ln b_{k-1} +1\over \delta -k + 1} n.
$$
\end{cor}

In particular, we obtain an upper bound for the total domination
number for any graph $G$ with $\delta\ge 1$: $
\gamma_{\mathrm t}(G) \le  {\ln\delta + 1 \over \delta} n.$

Finally, the case when $l=k$ in parametric domination is the $k$-tuple domination number $\gamma_{\times k}(G)$. 
Theorem \ref{main2} implies the following result, where 
$\tilde{b}_{k-1} = {b}_{k-1}+{b}_{k-2}= \pmatrix{\delta+1 \cr {k-1}}$.
These bounds can also be obtained from Theorem \ref{ub1} and Corollary \ref{cor1} if we put $r_i=1$ and $s_i=k$ for all $i=1,...,n$, in which case $r=1$, $k\le \delta$, $\theta=\delta-k+1$ and $B_{k-1}=\tilde{b}_{k-1}$.

\begin{cor} [Gagarin et al. \cite{gag1}]
For any graph $G$ with ${\widehat\delta}=\delta-k+1$ and $\delta\ge k$,
$$
\gamma_{\times k}(G) \le \left(1 - {{\widehat\delta} \over
({1+\widehat\delta})^{1+1/{\widehat\delta}} \; \tilde{b}_{k-1}^{1/{\widehat\delta}}}\right) n
$$
and
$$
\gamma_{\times k}(G) \le  {\ln(\delta-k+2) + \ln \tilde{b}_{k-1} +1 \over \delta-k+2} n.
$$
\end{cor}

\section{Threshold functions for multiple domination}

The bounds for multiple domination can be improved if we impose additional restrictions on graph parameters, i.e. by considering smaller graph classes. Such restrictions are called {\em threshold functions}. Caro and Roditty \cite{car1} and Stracke and Volkmann \cite{str1} were first who considered a threshold function for $k$-domination in the form $\delta \ge 2k-1$. 
For a slightly stronger threshold function Rautenbach and Volkmann \cite{rau1} 
found an interesting upper bound for the $k$-tuple domination number:

\begin{thm} [Rautenbach and Volkmann \cite{rau1}]
\label{rv} If $\delta\ge 2k\ln(\delta+1) -1$, then
$$
\gamma_{\times k}(G) \le   \left({k\ln(\delta +1)\over \delta + 1} +
\sum_{i=0}^{k-1}{k-i\over i!(\delta+1)^{k-i}}\right) n.
$$
\end{thm}

In the next theorem we consider a threshold function in the form $\delta\ge ck-1$, where $c>1$ is a constant. Although $c$ is not restricted from above, for given $k$ and $\delta$ the constant $c$ should not be taken as large as possible. The best approach would be to optimise $c$ for given $k$ and $\delta$ in such a way that the bound (\ref{const}) is minimised while $\delta\ge ck-1$ holds. We will deal with this later.
%and see that typically the constant $c$ is small. 

\begin{thm} \label{new1}
For any graph $G$ with $\delta\ge ck-1$, where $c>1$ is a constant, 
\begin{equation} \label{const}
\gamma_{\times k}(G) < \left({c \over \delta+1} + {1 \over e^{0.5k(c+1/c-2)}}\right) kn.
\end{equation}
\end{thm}

Let us consider a particular case of Theorem \ref{new1} when $c=3$, and compare it to Theorem \ref{rv} for graphs with $\delta\ge 20$.
We have
\begin{equation} \label{c=3}
\gamma_{\times k}(G) < \left({3 \over \delta+1} + {1 \over e^{2k/3}}\right) kn
\end{equation} 
for any graph $G$ with $\delta\ge 3k-1$. This bound is better than the bound of Theorem \ref{rv} if the former is less than the first term of the latter, i.e.
$$
\left({3 \over \delta+1} + {1 \over e^{2k/3}}\right) k < {k\ln(\delta+1) \over \delta+1},
$$
which is equivalent to 
$$
k > 1.5 \ln(\delta+1) -1.5 \ln[\ln(\delta+1)-3] = 1.5 \ln(\delta+1) (1-o(1)).
$$
Since Theorem \ref{rv} is applicable for $k\le (\delta+1)/(2\ln(\delta+1))$, we conclude that 
(\ref{c=3}) provides a better upper bound than Theorem \ref{rv} if
$$
1.5 \ln(\delta+1) (1-o(1))  < k \le  {\delta+1 \over 2\ln(\delta+1)},
$$
which is the largest part of the applicable interval. 

For example, if $\delta(G)=1000$, then 
Theorem \ref{rv} is applicable for $k\le 72$, whereas (\ref{c=3}) is applicable for $k\le 333$. 
Since $1.5 \ln(1001) -1.5 \ln[\ln(1001)-3]=8.3$, we obtain that the bound (\ref{c=3}) is stronger than the bound of Theorem \ref{rv} if
$$
9 \le k \le 72.
$$
If $k \le 8$, then Theorem \ref{rv} provides a better upper bound than (\ref{c=3}). However, we can try to optimise the constant $c$ in Theorem \ref{new1} for given $\delta$ and $k$ as follows. 
The right-hand side of the bound (\ref{const}) is minimised for $c$ satisfying the following equation:
$$
0.5k(c+{1\over c} -2) - \ln(0.5k(\delta+1))  =  \ln (1-{1\over c^2}).
$$
Now, replacing $\ln (1-{1\over c^2})$ by $-1/c^2$, we obtain the following cubic equation:
$$
kc^3 -2\left(k+\ln[0.5k(\delta+1)]\right) c^2 +kc +2 =0.
$$ 
The real root $c>1$ of this equation, which satisfies the condition $\delta \ge ck-1$, can be used in Theorem \ref{new1}. For example, if $k=5$ and $\delta=1000$, then the above cubic equation becomes
$$
c^3 - 5.13 c^2 +c+0.4=0.
$$
For this equation, the largest real root is $c=4.910$ (3 dp). Using this value of $c$ in Theorem \ref{new1}, we obtain $\gamma_{\times 5}(G) < 0.027 n$, whereas Theorem \ref{rv} produces the bound $\gamma_{\times 5}(G) < 0.035 n$.

It is not difficult to generalise Theorem \ref{new1} for parametric domination and $\left\langle \mathbf{r},\mathbf{s} \right\rangle$-domination. Let us denote 
$$\mu=\max\{k,l\}.$$

\begin{thm} \label{new2}
For any graph $G$ with $\delta\ge c\mu-1$, where $c>1$ is a constant, 
\begin{equation} 
\gamma_{k,l}(G) < \left({c \over \delta+1} + {1 \over e^{0.5\mu(c+1/c-2)}}\right) \mu n.
\end{equation}
\end{thm}

%It may be pointed out that the upper bound $\left\langle \mathbf{r},\mathbf{s} \right\rangle$-domination in the next theorem 
%is not dependent on $r$. 

\begin{thm} \label{new3}
For any graph $G$ with $(\delta+1)\tau\ge cs$, where $c>1$ is a constant, 
\begin{equation} 
\gamma {\left\langle \mathbf{r},\mathbf{s} \right\rangle}(G) < \left({c \over \delta+1} + {1 \over e^{0.5s(c+1/c-2)}}\right) sn.
\end{equation}
\end{thm}

Caro and Yuster \cite{car2} proved an important asymptotic result that if $\delta$ is much larger 
than $k$, then the upper bound for the total $k$-domination number is `close' to the bound
of Theorem \ref{arn}. More precisely, they proved the following:

\begin{thm} [Caro and Yuster \cite{car2}] \label{CY}
If $k<\sqrt{\ln\delta}$, then
\begin{equation}
\gamma_k^{\mathrm t}(G) \le {\ln \delta \over \delta}n (1+o_\delta(1)).
\end{equation}
\end{thm}

The same upper bound is therefore true for the $k$-tuple domination and $k$-domination numbers. 
The threshold function $k<\sqrt{\ln\delta}$ in Theorem \ref{CY} is indeed very strong,
but the corresponding bound is similar to that of Theorem \ref{arn}, which is best possible
in the class of all graphs. 
Let us consider 
a weaker but similar threshold function $k \le (1-c)\ln\delta$, where $0<c<1$ is a constant. 
The following explicit and asymptotic bounds are obtained: 

\begin{thm} \label{ln1}
For any graph $G$ with $k \le (1-c)\ln\delta$, where $0<c<1$ is a constant, 
\begin{equation}
\gamma_{\times k}(G) < 
\left({\ln \delta \over \delta+1} + {k \over \delta^{0.5c^2}}\right) n
\le {(1-c) \ln \delta \over \delta^{0.5c^2}} n (1+o_\delta(1)).
\end{equation}
\end{thm}

Theorem \ref{ln1} can be generalised for parametric domination and $\left\langle \mathbf{r},\mathbf{s} \right\rangle$-domination as follows.

\begin{thm} \label{ln2}
For any graph $G$ with $\mu \le (1-c)\ln\delta$, where $0<c<1$ is a constant, 
\begin{equation}
\gamma_{k,l}(G) < 
\left({\ln \delta \over \delta+1} + {\mu \over \delta^{0.5c^2}}\right) n.
\end{equation}
\end{thm}

Similar to Theorem \ref{new3} the upper bound in the next result does not depend on $\tau$. Note also that
$s\le (\delta+1) \tau$ holds because $s \le (1-c)\ln\delta < \ln\delta \le (\delta+1) \tau$, i.e. 
$\gamma{\left\langle \mathbf{r},\mathbf{s} \right\rangle}(G)$ is well defined.

\begin{thm} \label{ln3}
For any graph $G$ with $s \le (1-c)\ln\delta$, where $0<c<1$ is a constant, 
\begin{equation}
\gamma{\left\langle \mathbf{r},\mathbf{s} \right\rangle}(G) < 
\left({\ln \delta \over \delta+1} + {s \over \delta^{0.5c^2}}\right) n.
\end{equation}
\end{thm}

\end{document}